\documentclass{birkjour}

\usepackage{amssymb}
\usepackage{graphicx}
\usepackage{color}
\usepackage{enumitem}   

\newtheorem{Theorem}{Theorem}

\newtheorem{Lemma}[Theorem]{Lemma}
\newtheorem{Proposition}[Theorem]{Proposition}
\newtheorem{Corollary}[Theorem]{Corollary}

\newcommand{\R}{\mathbb R}
\newcommand{\Z}{\mathbb Z}
\newcommand{\E}{\mathbf{E}}
\newcommand{\pr}{\mathbf{P}}
\newcommand{\I}{\mathbf{I}}

\newcommand{\Eta}{\eta}
\newcommand{\erf}{\mbox{erf}}
\newcommand{\erfc}{\mbox{erfc}}
\newcommand{\K}{\mathbf{K}}
\newcommand{\Pf}{\mbox{Pf }}

\begin{document}

\title[Annihilating and Coalescing Random Walks on $\Z$]{Examples of Interacting Particle Systems on $\Z$ as Pfaffian Point Processes: Annihilating and Coalescing Random Walks}


\author[Garrod]{Barnaby Garrod}
\address{%
Mathematics Institute\\
University of Warwick\\
Coventry, CV4 7AL\\
UK}

\email{b.j.garrod@warwick.ac.uk}
\author[Poplavskyi]{Mihail Poplavskyi}
\address{%
Mathematics Institute\\
University of Warwick\\
Coventry, CV4 7AL\\
UK}

\email{m.poplavskyi@warwick.ac.uk}
\author[Tribe]{Roger P. Tribe}
\address{%
Mathematics Institute\\
University of Warwick\\
Coventry, CV4 7AL\\
UK}

\email{r.p.tribe@warwick.ac.uk}
\author[Zaboronski]{Oleg V. Zaboronski}
\address{%
Mathematics Institute\\
University of Warwick\\
Coventry, CV4 7AL\\
UK}

\email{olegz@maths.warwick.ac.uk}
\begin{abstract}
A class of interacting  particle systems on $\Z$, involving instantaneously annihilating or 
coalescing nearest neighbour random walks, are shown to be Pfaffian point processes 
for all deterministic initial conditions. As diffusion limits, explicit Pfaffian kernels are
derived for a variety of coalescing and annihilating Brownian systems. 
For Brownian motions on $\R$, depending on the initial conditions, the corresponding kernels are closely related to 
the bulk and edge scaling limits of the Pfaffian point process for real eigenvalues for the real
Ginibre ensemble of random matrices. For Brownian motions on $\R_{+}$ with absorbing or reflected boundary
conditions at zero new interesting Pfaffian kernels appear. We illustrate the utility of the Pfaffian 
structure by determining the extreme statistics of the right-most particle for the purely annihilating Brownian
motions, and also computing the probability of overcrowded regions for all models.
\end{abstract}
\maketitle
\section{Introduction and statement of key results}

Pfaffian point processes
arise in a number of contexts, for example the positions of eigenvalues of certain random matrix ensembles
(the Gaussian orthogonal ensemble, Gaussian symplectic ensemble, the real Ginibre ensemble), and in certain random 
combinatorial structures such as random tilings. 
In \cite{TZ:11}, the two systems of instantly coalescing, or instantly annihilating, Brownian motions, under a
maximal entrance law, were shown  to be Pfaffian point processes at any fixed time $t>0$. The aim of this paper
is to generalise this result in a number of ways:
\begin{enumerate}[label=(\roman*)]
 \item analogous particle systems on $\Z$; 
 \item mixed coalescent and annihilating systems;
 \item spatially inhomogeneous nearest neighbour motion; 
 \item general deterministic initial conditions. 
\end{enumerate}
The Pfaffian point process structure survives all of these changes. The key tool in the proof is a Markov duality. 

The introduction is organised as follows:  the discrete models are formulated in section \ref{s1.1}; the Pfaffian property is
stated in section \ref{s1.2}; the kernels for various diffusion limits are stated in section \ref{s1.3}; applications
of the Pfaffian property, and relation to other work are discussed in section \ref{s1.4}. 

\subsection{Mixed models.} \label{s1.1}
Interaction occurs when one particle jumps onto an already occupied site. 
 In a coalescent system there would be 
an instantaneous coalescence where the two particles merge to leave a single particle; in an annihilating system there 
would be an instantaneous annihilation where both particles disappear. We consider the following mixed system,
whose dynamics are informally described as follows. 

\noindent
\textbf{Particle Rules.} 
Between interactions all particles jump independently following a nearest neighbour random walk on $\Z$, jumping 
\[
\mbox{$x \to x-1$ at rate $q_x$, \hspace{.2in} and  \hspace{.2in}  $x-1 \to x$ at rate $p_x$.}
\]
The parameter $\theta\in [0,1]$ is fixed throughout, and when two particles interact they instantaneously annihilate with probability $\theta$ or coalesce with
probability $1-\theta$.  

We denote our particle systems as processes $(\Eta_t:t \geq 0)$ with values $\{0,1\}^{\Z}$, so that
$\Eta_t(x)=1$ indicates the presence of a particle at time $t$ at position $x$. 
The Markovian dynamics are encoded carefully in the generator (see (\ref{eq:processgenerator})).
We take the left and right jump rates $(p_x,q_x:x \in \Z)$ to be uniformly bounded, and then this generator determines the law of a unique Markov process, for any given initial condition in $\eta \in \{0,1\}^{\Z}$. We denote its law 
by $\pr_{\eta}$ and $\E_{\eta}$ on path space with canonical variables $(\Eta_t:t \geq 0)$. 


By choosing $\theta=0$ or $\theta=1$, our results apply to both purely coalescing and purely annihilating systems.
We know of three motivations for directly studying 
the mixed models $\theta \in (0,1)$.
\begin{itemize}
\item \textbf{Excitons in polymer chains.}
The kinetics of excitons (localised electronic excitations) along polymer chains are observed to exhibit 
sometimes coalescing collisions and sometimes annihilating collisions. Different polymer materials 
lead to different values of the parameter $\theta$ - see Henkel \cite{henkel} where three values of $\theta$ are 
listed for three different polymers. 
\item \textbf{Multi-valued voter models.} In the multi-valued voter model on $\Z$ (also known as the Potts model), started from one 
specific initial 
condition, the mixed systems arise as dual processes. 
Indeed, starting from product measure, where each of $q$ colours has equal
$1/q$ chance,  the domain walls that separate regions of sites with the same colour behave precisely as the 
mixed model above with the choice $\theta=1/(q-1)$.  This connection is explored in the amazing results of
Derrida et al. \cite{derridaetal1}, \cite{derridaetal2}  and we discuss these further in section \ref{s1.4}. 
\item \textbf{Coalescent mass models.}
Coalescent mass models are widely studied (see \cite{crtz}). The simplest model might be where massive particles perform simple coalescent random walks on $\Z$, and where masses add upon coalescence. Consider the case where
the initial masses are independent at different sites, and chosen uniformly from the set $\{0,1,2,\ldots,q-1\}$ for some fixed
(integer) $q \geq 1$.
Then the mass of particles in the future, modulo  $q$, remains uniformly distributed over $\{0,1,2,\ldots,q-1\}$. Moreover the
particles whose mass is exactly $r$ modulo $q$ perform a mixed coalescing/annihilating system with $\theta = 1/q$. This gives only a very partial description, and a full distribution of masses at multiple space points remains
an interesting problem. 
\end{itemize}

The spatially inhomogeneous version of the process here occurs in studies on reaction diffusion models with quenched disorder (see \cite{pierre}). The fact that the Pfaffian property extends to inhomogeneous nearest neighbour jump rates 
supports the suspicion discussed in {\cite{TZ:11} that in the continuum the free motion can be any 
continuous Markov process. 
\subsection{Statement of main theorem} \label{s1.2}
We have decided to recall the definition of a Pfaffian point process on $\Z$.
A variable $\Eta$ with values in $\{0,1\}^{\Z}$ can be thought of as a simple point process on $\Z$,
with $\Eta(x)=1$ corresponding to a particle at $x$.
In this context, the definition of a Pfaffian point process (see \cite{ppp}) is as follows: there exists a
matrix kernel $\K: \Z^2 \to \R^{2 \times 2}$, so that for any $n \geq 1$
\[
\E \left[ \Eta(x_1) \ldots \Eta(x_n) \right] = \Pf (\K(x_i,x_j): i,j \leq n) \quad \mbox{for any distinct $x_1, \ldots,x_n$.}
\]
$\Pf(A)$ is the Pfaffian of an antisymmetric matrix, and the   
appendix contains definitions and the simple properties of Pfaffians that we use.
The matrix kernel $\K$ can be written as
\[
\K(y,z) = \left( \begin{array}{cc} 
\K_{11}(y,z) & \K_{12}(y,z) \\
\K_{21}(y,z) & \K_{22}(y,z)
\end{array} \right)
\]
for $\K_{ij}:\Z^2 \to \R$, and these entries must satisfy the symmetry conditions
\[
\K_{ij}(y,z) = - \K_{ji}(z,y) \quad \mbox{for all $i,j \in \{1,2\}$ and $y,z \in \Z$}
\]
which ensure that the $2n \times 2n$ matrix $(\K(x_i,x_j): i,j \leq n)$, built out of the $n^2$ two-by-two blocks, is antisymmetric. 

The matrix kernel for our interacting particle system is constructed from a single scalar function $(K_t(y,z):y \leq z)$
defined as follows.  For $y,z \in \Z$ with $y\leq z$, and for $\eta \in \{0,1\}^{\Z}$, we define
\[
\eta[y,z) = \sum_{y \leq x<z} \eta(x) \quad \mbox{if $y<z$,}
\]
and $\eta[y,y)=0$, and we define the 'spin pair' by
\[
\sigma_{y,z}(\eta) = (-\theta)^{\eta[y,z)}.
\]
We use the convention that $0^0=1$ so that when $\theta=0$ the spin reduces to the  indicator of an empty interval, that is
$\sigma_{y,z}(\eta) = \I(\eta[y,z)=0)$. We now set
\begin{equation} \label{eq:scalarkernel}
K_t(y,z) = \E_{\eta} \left[ \sigma_{y,z}(\Eta_t) \right], \quad \mbox{for $t \geq 0$, $y,z \in \Z$ with $y \leq z$.}
\end{equation}
We also need the difference operators $D^+$ and $D^-$, defined for $f:\Z \to \R$ by
\[
D^+ f(x) = f(x+1)-f(x), \qquad D^- f(x) = f(x-1)-f(x).
\]
\begin{Theorem} \label{T1}
For any initial condition $\eta \in \{0,1\}^{\Z}$, and  at any fixed time $t \geq 0$, the variable $\Eta_t$ is a Pfaffian point process on $\Z$ with 
kernel $\K$ given, for $y<z$, by
\begin{equation} \label{discretekernel}
\K(y,z) = \frac{1}{1+\theta} \left( \begin{array}{cc} 
K_t(y,z) & - D^+_2K_t(y,z) \\
- D^+_1K_t(y,z) & D^+_1 D^+_2 K_t(y,z)
\end{array} \right),
\end{equation}
and $\K_{12}(y,y) =  \frac{-1}{1+\theta} \; D^+_2 K_t(y,z)\mid_{z=y}$, and other entries determined by the symmetry conditions.
(The notation $D^{\pm}_i$ means that the operator $D^{\pm}$ is applied in the $i$th variable.)  
\end{Theorem}

\textbf{Remarks}. 

\textbf{1.} \textbf{Random initial conditions.} 
For random initial conditions the law of $\eta_t$ is not in general a Pfaffian point process, though by conditioning
on the initial condition $\eta_0$ the correlation functions can always be expressed as the expectation of a 
Pfaffian with a random kernel $\K^{(\eta_0)}$ depending on $\eta_0$:
\begin{equation} \label{rii}
\rho^{(n)}_t(x_1,\ldots,x_n) = \E \left( \Pf \left(\K^{(\eta_0)}(x_i,x_j): i,j \leq n \right) \right).
\end{equation}
For certain random initial conditions, including the natural case when the sites $(\eta_0(x):x \in \Z)$ are independent,
the expectation can be taken inside the Pfaffian and the process does remain a Pfaffian point process
See the first remark at the end of section \ref{s2}.

The simplest random initial condition is Bernoulli-$1/2$. In this case, the kernel can be written explicitly by solving
a linear system of ODE's as explained in Section \ref{s2}. The answer is (\ref{discretekernel}), where
\begin{eqnarray}
K_t(x,y)=2+2e^{-2t}\sum_{l=1}^{\infty} \left(I_{y-x+l}(2t)-I_{y-x-l}(2t)\right)\mbox{ for } x\leq y, ~t>0,
\end{eqnarray}
where $I_n$ is the Bessel function of the imaginary argument defined via
\[
e^{\frac{1}{2}x(\lambda+\lambda^{-1})}=\sum_{n \in \Z} \lambda^n I_n(x).
\]
This answer dates back to the seminal $1963$ paper by Glauber \cite{glauber}, where the kinetic
Ising chain was introduced and analysed. Indeed, the duality function for the annihilating case is just
a two-point spin-spin function computed explicitly in Glauber's paper.

\textbf{2.} \textbf{Thinning.} The parameter $\theta$ enters into the kernel only as a scalar
mulitiplier.  
Instantly coalescent systems and instantly annihilating systems are related by a 
well known thinning relation (see \cite{TZ:11} section 2.1 or \cite{thinning}).  This extends to our mixed systems
as follows: consider a two colour system of particles, red $R$ and blue $B$, that move independently between
reactions, and at reaction times (when one particle lands on top of another) transform via the rules
\[
R+R \to R, \quad R+B \to B, \quad   B+B \to \left\{ 
\begin{array}{cl}
R & \mbox{with probability $\theta$,} \\
B & \mbox{with probability $1-\theta$.} 
\end{array} \right.
\] 
Note that the full system of particles, where one ignores colours, is a 
coalescing system, but the blue particles alone follow the mixed model we study. However, if initially particles 
are coloured 
\[
\mbox{blue with probability $\frac{1}{1+\theta}$ and red with probability $\frac{\theta}{1+\theta}$.}
\]
then this property is preserved at all subsequent times. Indeed one can check that each single collision preserves
property. Thus at any fixed time the mixed system is a thinning, by the factor $1/(1+\theta)$ of the full coalescing system. 
Thinning also acts naturally on Pfaffian point processes, changing the underlying kernel by the same factor.  
However this connection  seems to relate the two systems only when the initial conditions are 
similarly related by thinning, and so does not apply to the deterministic initial conditions stated in the Theorem.

\textbf{3.} \textbf{Strong thinning.} 
The Pfaffian point processes with kernel (\ref{discretekernel}) but where the scalar
multiplier $1/(1+\theta)$ is replaced by an arbitrary $\lambda \in (0,1)$ also arise. When $\lambda \in (0,\frac12)$ these 
correspond to thinning a coalescing system by more than the factor needed to reach  annihilating systems. 
Consider a two colour system as described above, but with reactions
\[
R+R \to R, \quad B+B \to B, \quad   B+R \to \left\{ 
\begin{array}{cl}
R & \mbox{with probability $\frac12$,} \\
B & \mbox{with probability $\frac12$.} 
\end{array} \right.
\] 
If we initially choose colours independently as
\[
\mbox{blue with probability $\lambda$ and red with probability $1-\lambda$}
\]
then this property is preserved at all later times. Thus the sub-population of blue particles is a thinning, by the
factor $\lambda$, of a coalescing system and hence Pfaffian with the coalescing kernel thinned by the factor
$\lambda$. 
\subsection{Continuous limits} \label{s1.3}
The scalar function $K_t(y,z)$ that underlies the Pfaffian matrix kernel $\K$ can be characterised as the solution to
a system of differential equations indexed over part of the lattice. We define a one-particle generator $L$, acting on $f:\Z \to \R$, by 
\begin{equation}
Lf (x)  = q_x D^+f(x) + p_{x} D^-f(x).   \label{eq:particlegenerator}
\end{equation}
The intuition is that $L$ is the generator for a single dual particle.
We will show that the function $(K_t(y,z):t \geq0, \, y,z \in \Z, \,  y < z)$ is the unique bounded solution to the equation
\begin{equation} \label{eq:ode2}
\left\{ \begin{array}{rcll} 
\partial_t K_t(y,z) &=& (L_y + L_z) K_t(y,z) & \mbox{for $y<z, \, t >0$,} \\
K_t(y,y) & = & 1 & \mbox{for all $y, \, t >0$,}  \\
K_0(y,z) & = & \sigma_{y,z}(\eta) & \mbox{for $y \leq z$.} 
\end{array} \right.
\end{equation}
(The notation $L_y$ is used to indicate that the operator $L$ acts on the $y$ variable.)

This differential equation characterisation (\ref{eq:ode2}) lends itself naturally to asymptotic analysis, where 
its large time and space behaviour is determined by a similar limiting continuum kernel $K^{(c)}_t(y,z)$ solving a continuum PDE. We can in several natural cases solve these limiting continuum equations explicitly, and therefore add to the growing zoo of concrete known kernels for Pfaffian point process on $\R$. 

In the cases below we use the diffusive scaling 
\begin{equation} \label{scaledprocess}
X^{(\epsilon)}_t = \sum_{x \in \Z}  \eta_{\epsilon^{-2}t}(x) \delta_{\epsilon x}
\end{equation}
and we check, at a fixed $t$, that $X^{(\epsilon)}_t \to X_t$ as $\epsilon \to 0$ (considered as random locally finite point measures on $\R$ with the topology of vague convergence). By choosing $p_x,q_x,\theta$ and the initial conditions
we establish the Pfaffian property for various continuum systems at fixed times. In each case the limiting point measure
$X_t$ is a Pfaffian point process on with a kernel of the form
\begin{equation} \label{ctmkernel}
\K^{(c)}(y,z) = \frac{1}{1+\theta} \left( \begin{array}{cc} 
K^{(c)}_t(y,z) & -\partial_z K^{(c)}_t(y,z) \\
-\partial_y K^{(c)}_t(y,z) & \partial^2_{yz} K^{(c)}_t(y,z)
\end{array} \right), \quad \mbox{for $y<z$,} 
\end{equation}
and $\K^{(c)}_{12}(y,y) =  \frac{-1}{1+\theta} \; \partial_2 K^{(c)}_t(y,z)\mid_{z=y}. $

We record some specially chosen cases where
$K^{(c)}_t(y,z)$ can be found explicitly. The limits can also be identified as the law at time $t$ for a suitable system of 
reacting Brownian particles. Particularly simple kernels appear for the Poisson 
initial distribution of particles in the limit of infinite intensity. We refer to such a limit
as the maximal entrance law, see 
section \ref{s3} for a more formal discussion of entrance laws for mixed systems. We emphasise that the Pfaffian property holds for all deterministic initial conditions - but the
maximal initial condition leads to a simple initial condition of the PDE's determining the kernel and hence simple explicit solution formulae. More importantly, as explained in \cite{TZ:11} for purely coalescing or annihilating systems, the maximal entrance law is distinguished in that 
the solutions are then invariant under diffusive rescaling. Moreover the solutions from a large class of initial 
distributions becomes attracted under diffusive rescaling to the solution started form the the maximal entrance law. 
\begin{Theorem} \label{T2} 
Fix $t>0$ throughout. 
Recall $\erfc(x) = 1- \erf(x)$ for $\erf(x) = \frac{2}{\sqrt{\pi}}\int^x_0 \exp(- u^2) du$. 
\begin{enumerate}
\item[(A)] \textbf{Brownian motions on $\R$ - maximal entrance law.}

Take $p_x=q_x=1$ for all $x$ and the initial condition $\eta_0 \equiv 1$. The diffusion limit $X_t$ has the Pfaffian kernel (\ref{ctmkernel}) on $\R$ with
\begin{equation}  \label{kernelbulk}
K^{(c)}_t(y,z) = \erfc \left(\frac{z-y}{2 \sqrt{2t}}\right).
\end{equation}
Moreover, this is the kernel for mixed coalescing/annihilating Brownian motions on $\R$ at time $t$ under the maximal entrance law. 
\item[(B)] \textbf{Brownian motions on $\R$ - half-space maximal entrance law.}

Take $p_x=q_x=1$ for all $x$ and half-space initial conditions: $\eta_0(x) = 1$ for $x \leq 0$ and $\eta_0(x) = 0$ for $x>0$.
The diffusion limit $X_t$ has the Pfaffian kernel (\ref{ctmkernel}) on $\R$ with
\begin{equation}  \label{kerneledge}
K^{(c)}_{t}(y,z) = 1 + \int^{\frac{z}{2\sqrt{t}}}_{\frac{y}{2\sqrt{t}}} \int^{\frac{y}{2\sqrt{t}}}_{-\infty} \frac{(u-v)}{\sqrt{2 \pi}} e^{-\frac{(u-v)^2}{2}} \erfc\left(\frac{u+v}{\sqrt{2}}\right)  du dv.
\end{equation}
Moreover, this is the kernel for mixed coalescing/annihilating Brownian motions on $\R$ under the half-space 
maximal entrance law. 
\item[(C)] \textbf{Killed Brownian motions on $[0,\infty)$ - maximal entrance law.}

Take
\[
q_x = \left\{ \begin{array}{ll}
1 & \mbox{for $x \geq 1$,} \\ 
2 & \mbox{for $x=0$,} \\
0 & \mbox{for $x <1$,}
\end{array} \right.
 \qquad 
p_x = \left\{ \begin{array}{ll}
1 & \mbox{for $x \geq 1$,} \\ 
0 & \mbox{for $x <1$.}
\end{array} \right.
\]
Also take $\Eta_0(x) = 0$ for $x \leq 0$ and $\Eta_0(x) = 1$ for $x \geq 1$. 
The diffusion limit $X_t$ has the Pfaffian kernel (\ref{ctmkernel}) on $(0,\infty)$ with
\begin{equation}  \label{kernelkilled}
K^{(c)}_t(y,z) = 1 - \erf\left(\frac{z+y}{2\sqrt{2t}}\right) \erf\left(\frac{z-y}{2\sqrt{2t}}\right).
\end{equation}
Moreover, this is the kernel for mixed coalescing/annihilating Brownian motions on $[0,\infty)$, killed at $\{0\}$, under the half-space maximal entrance law. 
\item[(D)]  \textbf{Reflected motions on $[0,\infty)$ - maximal entrance law.}

Take 
\[
q_x = \left\{ \begin{array}{ll}
1 & \mbox{for $x \geq 1$,} \\ 
0 & \mbox{for $x <1$,}
\end{array} \right.
 \qquad 
p_x = \left\{ \begin{array}{ll}
1 & \mbox{for $x \geq 1$,} \\ 
0 & \mbox{for $x <1$.}
\end{array} \right.
\]
Also take $\Eta_0(x) = 1$ for $x \geq 0$ and $\Eta_0(x) = 0$ for $x <0$,
The diffusion limit $X_t$ has the Pfaffian kernel (\ref{ctmkernel}) on $[0,\infty)$ with $K^{(c)}_{t}(y,z)$ given by
\begin{equation} \label{kernelreflected}
1 
-\frac14  \int^{\frac{z}{2\sqrt{t}}}_{\frac{y}{2\sqrt{t}}} du \int^{\infty}_{\frac{y}{2\sqrt{t}}} dv \;
\left(
 \erf^{\,\,''} \left( \frac{u-v}{\sqrt{2}} \right)   \erf \left( \frac{u+v}{\sqrt{2}} \right) +
 \erf^{\,\,''} \left( \frac{u+v}{\sqrt{2}} \right)   \erf \left( \frac{u-v}{\sqrt{2}} \right) \right). 
\end{equation}
Moreover, this is the kernel for mixed coalescing/annihilating Brownian motions on $[0,\infty)$, reflected at $\{0\}$, 
under the half-space maximal entrance law. 
\end{enumerate}
\end{Theorem}
The final two examples illustrate the fact that we allow spatial inhomogeneities
in our discrete models. One point intensities $ \rho_t(y) = - \partial_z K_t^{(c)}(y,z) |_{z=y}$
can then be read off from the Pfaffian kernels yielding for the absorbing case
\[
\rho_t(y) = \frac{1}{1+\theta} \frac{1}{\sqrt{2\pi t}} \erf\left(\frac{y}{\sqrt{2t}}\right)
\]
and for the reflecting case
\[
\rho_t(y) = \frac{1}{1+\theta}  \left(\frac{1}{\sqrt{2 \pi t}} \erf \left( \frac{y}{\sqrt{2t}} \right) + \frac{1}{\sqrt{\pi t}}
e^{-\frac{y^2}{4t}} \erfc \left( \frac{y}{2 \sqrt{t}} \right) \right).
\]
The figure below illustrates these for the coalescing case $\theta=0$ at $t=1/4$. Note that both intensities
converge as $y \to \infty$ to that of example (A), the free Brownian motions. 
\begin{figure}[t!]
\centering
\includegraphics[width=0.55\textwidth]{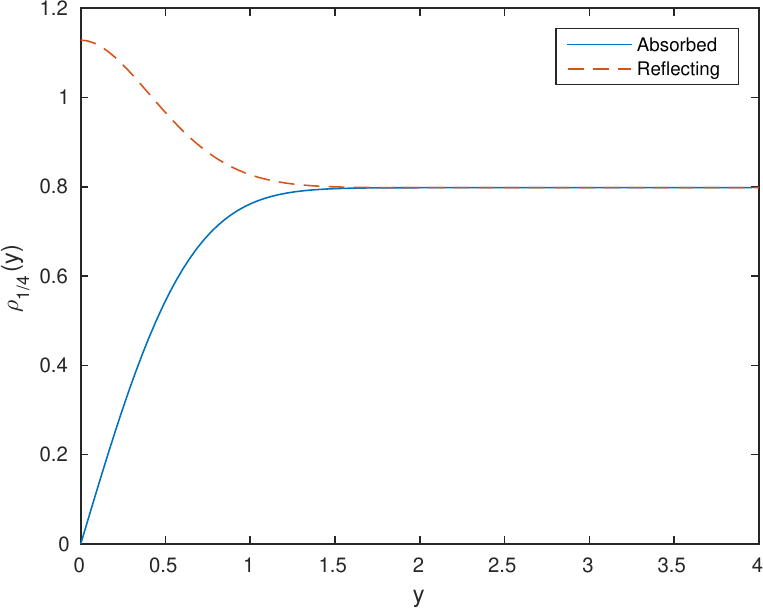}
\caption{Particle density at time $t=1/4$ for coalescing Brownian motions killed at the
origin (example (C)), and reflecting at the origin (example (D)). For the absorbing
case, particle density vanishes at the origin due to the trap at $0$.}
\end{figure}
\subsection{Related work and applications.} \label{s1.4}
\textbf{The real Ginibre random matrix ensemble.}
The real Ginibre random matrix ensemble is the $N \times N$ random matrix formed by taking
independent $N(0,1)$ Gaussian variables as entries. The positions of the eigenvalues form a 
Pfaffian point process (see \cite{BS:09}, \cite{FN:07}, \cite{SW:08}). 
Letting $N \to \infty$ there is a limiting kernel describing the
eigenvalues in the bulk of the spectrum. The paper \cite{TZ:11}, thanks to the comments of 
a referee, noted that the kernel for the
real eigenvalues in the bulk can be written in the form (\ref{ctmkernel}) where $K_t(y,z)$ is given by 
(\ref{kernelbulk}) when $\theta=1$, pure annihilation, and $t=1/4$. 

A similar large $N$ limit can be taken for the real eigenvalues in the real Ginibre ensemble near the right hand edge of the spectrum. 
Again the limiting kernel can be written in the form (\ref{ctmkernel}) where $K_t(y,z)$ is now given by 
(\ref{kerneledge}) when $\theta=1$ and $t=1/4$. 
(The original version of \cite{BS:09} has a slightly incorrect derivation of one of the four 
entries in the limiting edge kernel $\K^{Ginibre}_{Edge} (y,z)$. However
this can be easily corrected and all four entries agree with the kernel above (see the erratum \cite{erratum})). 
We state these observations as a corollary.

\begin{Corollary}
The large $N$ kernels for the real eigenvalues of the real Ginibre ensemble, in the bulk or at the edge, 
agree with the
kernels for annihilating Brownian motions, at time $t=1/4$, started from the maximal entrance law or the 
half-space maximal entrance law. 
\end{Corollary}
\noindent
\textbf{The probability of overcrowded regions.}
The Pfaffian structure for mixed models is well suited for the analysis of the following
natural question: what is the probability of finding a configuration of particles separated by
distances much smaller than the typical distance? 
To be more precise, consider the system
of coalescing-annihilating Brownian motions on $\R$ started from the maximal entrance law. 
The corresponding single-time distribution of particles is a Pfaffian point process  with the 
kernel given by (\ref{ctmkernel}) and (\ref{kernelbulk}). Particle density decays with time
as $t^{-1/2}$, meaning that the typical inter-particle separation is $O(\sqrt{t})$.
What is the probability density for finding particles in positions $x_1, x_2, \ldots, x_n$
such that 
\[
\max_{1\leq i,j \leq n}|x_i-x_j|=L<<\sqrt{t}\mbox{ ?}
\]
The answer is a simple corollary
of the corresponding answer for purely coalescing systems obtained in \cite{TZ:11}. 
Using Theorem $1$ of the above paper, (\ref{ctmkernel}) and the fact that for pure coalescence $\theta=0$, we conclude
that
\begin{eqnarray}\label{eqn:dens_dens}
\rho_t(x_1, x_2, \ldots, x_n)=\frac{C_n}{(1+\theta)^n}t^{-\frac{n(n+1)}{4}}
\mid \Delta(x_1, x_2,\ldots, x_n) \mid \left(1+O\left(\frac{L}{\sqrt{t}}\right)\right),
\end{eqnarray} 
where $\Delta(x_1, \ldots, x_n)$ is the Vandermonde determinant and $C_n>0$ is a
universal (i.e. time and $\theta$-independent) constant computed in \cite{TZ:11}. 
Integrating (\ref{eqn:dens_dens}) over the box $[0,L]^n$ and using the fact that $\rho_t$ is the density for the
factorial moments of the total number of particles $X_t[0,L]$ in the interval, we find that
\begin{eqnarray}\label{eqn:dens_prob}
P[ X_t[0,L]=n]=
\frac{M_n}{(1+\theta)^n}\left(\frac{L}{\sqrt{t}}\right)^{\frac{n(n+1)}{2}}
 \left(1+O\left(\frac{L}{\sqrt{t}}\right)\right),
\end{eqnarray}
where $M_n>0$ does not depend on $t,\theta$ and $L$. Formula (\ref{eqn:dens_prob}) 
quantifies the negative dependence between coalescing-annihilating Brownian motions: namely, we see that
\begin{eqnarray}\label{eqn:dens_ratio}
\frac{P[ X_t[0,L]=n]}
{\left(P[X_t[0,L]=1]\right)^n}
=N_n\left(\frac{L}{\sqrt{t}}\right)^{\frac{n(n-1)}{2}}
\left(1+O\left(\frac{L}{\sqrt{t}}\right)\right),
\end{eqnarray}
where $N_n>0$ does not depend on $t,\theta$ and $L$.
\\
\\
\noindent
\textbf{Gap probabilities.}
The opposite 
question often asked in the context of the determinantal and Pfaffian point processes is the problem of gap probabilities
\[
p(t,[L,R]) = P[ \mbox{there are no particles in an interval $[L,R]$ at time $t$}]
\]
and their asymptotics as $R-L \to \infty$. A particularly interesting special case of gap probability is $p(t, [L,\infty))$
under the half space initial condition, which gives the law of the right most particle. 

These probabilities can be 
expressed in terms of the Pfaffian kernels as Fredholm Pfaffians, and their study for the 
Pfaffian point processes arising from classical random matrix 
ensembles (GOE and GSE) is described in \cite{TW}. 
For the interacting particle systems described here, bulk gap probabilities were studied by 
Derrida and Zeitak \cite{derridaetal2}. They were motivated by the connection with the Potts model described above, 
and this led them to analyse the gap probability under a product Bernoulli initial condition, in the translation invariant case.  The leading term of their asymptotic might be expressed as 
\begin{equation} \label{derridaasymptotic}
\lim_{t \to \infty} p(t, t^{1/2}[L,R]) \approx e^{-A(\theta)(R-L)} \quad \mbox{for $R-L$ large}
\end{equation}
 where 
\[
A(\theta) = \frac14 (1+\theta) \sum_{n=1}^{\infty} \frac{1}{n^{3/2}} \left(\frac{4 \theta}{(1+\theta)^2}\right)^n.
\]
Although they do not identify a Pfaffian point process, the arguments in \cite{derridaetal2} exploit Pfaffians, and manipulations similar to those exploited for random matrix ensembles, and this paper was a motivation in our search for an underlying Pfaffian point process for these mixed models.  In a subsequent paper we will show that 
these techniques, starting from the Fredholm Pfaffian, will confirm the asymptotic 
(\ref{derridaasymptotic}) for all $\theta \in [0,1]$ for Bernoulli initial conditions, 
and then extend them to more general initial conditions.

One non-translation invariant case has already been investigated, namely the 
edge gap probability $p(t, [L,\infty))$ for the annihilating particle system. This corresponds to the
largest real eigenvalue in the real Ginibre ensemble and has been studied in \cite{RS} and \cite{PTZ}.
The main theorem in \cite{PTZ} can be expressed as follows.
\begin{Theorem}\label{thm4}
Consider the system of instantaneously annihilating Brownian motions on $\R$ with half-space
maximal initial conditions. Then 
\begin{equation}
p(t,[L,\infty)) = 
\left\{ \begin{array}{ll}
1 - \frac{1}{4}\erfc(\frac{L}{\sqrt{4t}}) + O(e^{-\frac{L^2}{4t}}) &  \mbox{for $\frac{L}{\sqrt{t}} \to \infty$,} \\
 \exp( \frac{1}{2 \sqrt{2 \pi}} \zeta(\frac32) \frac{L}{\sqrt{4t}} + O(1))
 & \mbox{for $\frac{L}{\sqrt{t}} \to- \infty$.}
 \end{array} \right.
\end{equation}
\end{Theorem}
Notice that the applicability condition for the above statement involves 
only the ratio $L/\sqrt{t}$, so the theorem still applies even if $t$ is allowed to scale with $L$.
The techniques used to prove Theorem \ref{thm4} start with the Pfaffian point process kernel, but exploit a representation for the 
terms in the Fredholm Pfaffian in terms of a discrete time simple random walk, reducing the asymptotic to 
a simpler asymptotic for a single random walk. We believe the same techniques will apply for the mixed 
coalescing/annihilating models. 

It is interesting to compare the above tail gap problem with the results 
connecting the statistics for interacting particles in the KPZ universality class, started with half-space
initial conditions, with the Tracy-Widom distribution corresponding to the behaviour of the largest eigenvalue
for Hermitian matrix models.
The fact that annihilating Brownian motions appear as the limit for a large number of discrete 
interacting particle models, suggests that the distribution of the largest real eigenvalue for the 
real Ginibre ensemble may define an interesting universality class associated with 
statistics for non-Hermitian matrices.  

The underlying Pfaffian structure should be useful in the study of many questions concerning these models,
for example multi-time correlations. 
The extended Pfaffian property found for annihilating systems in \cite{TYZ} was shown in \cite{thesis} to extend to
for the discrete space inhomogeneous process studied here, solely in the annihilating case. 
The persistence exponents found by Derrida, Bray, and Godreche \cite{derridaetal1} suggest that
Pfaffian structure will be useful for time dependent properties as well. 
\section{Proof of the main result.}  \label{s2}
We start with a summary of the main steps, which follow similar lines as \cite{TZ:11}. 
The key tool is a Markov duality. Indeed for any $n \geq 1$ the product of $n$ spin pairs
$\eta \to \sigma_{y_1,y_2}(\eta) \ldots \sigma_{y_{2n-1},y_{2n}}(\eta)$
is a suitable Markov duality function, as shown in Lemma \ref{lemma:generator}. Exploiting this allows us to calculate 
the expectations
\[
\E_{\eta} \left[ \sigma_{y_1,y_2}(\eta_t) \ldots \sigma_{y_{2n-1},y_{2n}}(\Eta_t) \right]
\]
as the solutions of $2n$ dimensional (spatially inhomogeneous) lattice heat equations. 
This is similar to the Markov dualities used in \cite{BCS} to study the ASEP and q-TASEP models. 
The dual process can be taken to be (a spatially inhomogeneous version of) the one dimensional 
Glauber spin chain (see the Remark after Lemma \ref{lemma:generator}). This model is known 
to be solvable by mapping to a system of free fermions operators (see \cite{glauberpfaffs}). Fermions are naturally 
associated to Pfaffians, and it turns out that the duality expectations are given by $2n \times 2n$ Pfaffians of a matrix built from a scalar kernel $K_t(y_i,y_j)$, as shown in Lemma \ref{lemma:scalarpfaffian}. 
The final step is to reconstruct the particle intensities from the product spin expectations. This is possible via the identity
\begin{equation} \label{eq:reconstruct}
\eta(y) = \frac{1-\sigma_{y,y+1}(\eta)}{1+\theta}
\end{equation}
This leads to a linear reconstruction formula for the $n$ point intensity in terms of the $2n \times 2n$ Pfaffians for the 
product of $n$ spin pairs. The Pfaffian structure is preserved by the reconstruction formula, but the matrix breaks into $2 \times 2$ blocks
corresponding to the spin pairs, and this yields the desired matrix kernel $\K(y_i,y_j)$. 

We now present the details. 
The generator of the process is given, for $F:\{0,1\}^{\Z}\to \R$ that depend on finitely many coordinates, by 
\begin{eqnarray}
 \mathcal{L}  F(\eta) & = & \sum_{x\in\Z} q_x \left (\theta F(\eta_{x,x-1}^a)+(1-\theta)F(\eta^c_{x,x-1})-F(\eta)  \right) \nonumber \\
&& \hspace{.2in}  +\sum_{x\in\Z} p_{x} \left(\theta F(\eta^a_{x-1,x})+(1-\theta)F(\eta^c_{x-1,x}) - F(\eta)\right), \label{eq:processgenerator}
\end{eqnarray}
where $\eta^a_{x,y}$ (respectively $\eta^c_{x,y}$) is the new configuration after a jump from site $x$ to $y$ followed by instantaneous annihilation
(respectively coalescence). These are defined, when $x \neq y$, by 
\[
\left\{
\begin{array}{l}
\eta^a_{x,y}(z) = \eta^c_{x,y}(z)= \eta(z) \quad \mbox{for $z \not \in \{x,y\}$,} \\ 
\eta^a_{x,y}(x) = \eta^c_{x,y}(x) = 0, \\
\eta^a_{x,y}(y) = (\eta(x) + \eta(y)) \; \mbox{mod}(2), \\
\eta^c_{x,y}(y) = \min\{1, \eta(x)+\eta(y)\}.
\end{array} \right.
\]

The key to the argument is the following Markov duality function.
For $n \geq 1$ and $y =(y_1,\ldots,y_{2n})$ with $y_1 \leq y_2 \leq  \ldots \leq y_{2n}$ we define the product spin function by 
\[
\Sigma_y(\eta) = \prod_{i=1}^n \sigma_{y_{2i-1},y_{2i}}(\eta). 
\]
Note that $\Sigma_y(\eta)$ depends only on finitely many coordinates of $\eta$ and so lies in the domain of the generator $\mathcal{L}$. 
The Markov duality is encoded in the following generator calculation.
\begin{Lemma} \label{lemma:generator}
 For $y_1<y_2<\dots<y_{2n}$ the action of the particle generator $\mathcal{L}$ on $\Sigma_{y}(\eta)$ is
 \[
  \mathcal{L} \, \Sigma_{y}(\eta)=\sum_{i=1}^{2n} \, L_{y_i} \Sigma_y(\eta),
 \]
 where $L_{y_i}$, given by (\ref{eq:particlegenerator}), acts on the coordinate $y_i$ in $\Sigma_y$.
\end{Lemma}

\textbf{Remark.}  We do not make use of a dual Markov process, but this lemma could be cast into the standard framework
(see Ethier and Kurtz \cite{EK:86} chapter 4) relating two Markov processes. The dual process can be taken to be a finite system 
of particles with motion generator $L$ that are instantly annihilating (with state space the disjoint union  $\cup_{m=0}^n \R^{2m}$). 
However this annihilating system describes the motion of domain walls in (a spatially inhomogeneous version of) 
the Glauber dynamics for the Ising spin chain \cite{glauber} and the dual process could also be taken to be this spin chain. 
The formulae connecting a set of spins $(\sigma(x) \in \{-1,+1\}, x \in \Z)$ to the positions a domain wall separating different spins at $x$ and $x+1$ 
are
\[
\eta(x) = \frac{1-\sigma(x) \sigma(x+1)}{2}, \qquad (-1)^{\eta[x,y)} = \sigma(x) \sigma(y).
\]
We do not exploit the link between the spin chain and annihilating systems but it is the origin of our use of the term 'spin-pair'. 
The fact that the dual process both for coalescing or annihilating systems can be taken to be the spin chain perhaps explains why both systems are Pfaffian. 

\textbf{Proof.}
Each term in the generator $\mathcal{L}$ involves a modified configuration, $\eta^a_{x,y}$ or $\eta^c_{x,y}$, which differs form $\eta$ on at most two 
neighbouring sites. The condition that $y_i < y_{i+1}$ ensures that this modified configuration will agree with $\eta$ on all but at most one 
of the intervals $[y_{2i-1},y_{2i})$, and hence the value of at most one of the 
spin pairs $\sigma_{y_{2i-1},y_{2i}}$ will change. This allows us to separate the action of the generator as follows
\begin{equation} \label{eq:generator_separate}
 \mathcal{L} \, \Sigma_y (\eta)=\sum_{i=1}^n \left(\prod_{j=1,j\neq i}^n \sigma_{y_{2j-1},y_{2j}}(\eta) \right) \mathcal{L} \, \sigma_{y_{2i-1},y_{2i}}(\eta).
\end{equation}
We turn our attention to a single spin pair $\sigma_{y,z}(\eta)$. Fix $y<z$ and consider the part of the generator
\[
 \sum_{x\in\Z} q_x \left (\theta \sigma_{y,z}(\eta^a_{x,x-1})+(1-\theta) \sigma_{y,z}(\eta^c_{x,x-1})- \sigma_{y,z}(\eta)  \right)
\]
corresponding to left jumps. The terms in this sum indexed by $x \leq y$ and by $x \geq z+1$ are zero, as the modified
configuration is unchanged in the interval $[y,z)$.
 The terms corresponding to $y+1 \leq x \leq z-1$ are also zero since
we claim that 
\[
\theta \sigma_{y,z}(\eta^a_{x,x-1})+(1-\theta) \sigma_{y,z}(\eta^c_{x,x-1})- \sigma_{y,z}(\eta) = 0.
\]
Indeed, since $\{x-1,x\} \in [y,z)$, the left hand side is proportional to
\[
\theta \, (-\theta)^{\eta(x-1)+\eta(x) \; \mbox{mod(2)}} + (1-\theta) \, (-\theta)^{\min(1,\eta(x-1)+\eta(x))}
- (-\theta)^{\eta(x-1)+\eta(x)}
\] 
and checking the three cases $\eta(x-1)+\eta(x) \in \{0,1,2\}$ shows that this is always zero.
(This identity is similar to a key quadratic identity (relating the exponential parameter $\tau$ to the 
ASEP parameters $p,q$) that lies behind the ASEP dualities 
in section 4.1 of \cite{BCS}.)
Thus jumps between sites both lying outside or both lying inside the interval $[y,z)$ give zero contribution to the generator and
only the terms when $x=y$ or $x=z$, where a jump crosses an endpoint of the interval, contribute. It was this key property that
we looked for when trying to find a duality function.

For the two remaining terms, when $x=y$ or $x=z$, two similar short calculations lead to 
\begin{eqnarray*}
 \theta \sigma_{y,z} (\eta^a_{y,y-1})+(1-\theta) \sigma_{y,z}(\eta^c_{y,y-1})- \sigma_{y,z}(\eta) &=& \sigma_{y+1,z}(\eta)- \sigma_{y,z}(\eta) = D^+_y \sigma_{y,z}(\eta),\\
 \theta \sigma_{y,z} (\eta^a_{z,z-1})+(1-\theta) \sigma_{y,z} (\eta^c_{z,z-1})- \sigma_{y,z}(\eta)&=&\sigma_{y,z+1}(\eta)-\sigma_{y,z} (\eta) = D^+_z \sigma_{y,z}(\eta).
\end{eqnarray*}
Collecting the terms from all possible $x$ gives 
\[
 \sum_{x\in\Z} q_x \left (\theta \sigma_{y,z}(\eta^a_{x-1,x})+(1-\theta) \sigma_{y,z}(\eta^c_{x-1,x})- \sigma_{y,z}(\eta)  \right)
= q_y D^+_y \sigma_{y,z}(\eta) + q_z D^+_z \sigma_{y,z}(\eta).
\]
Similar calculations for the terms corresponding to right jumps show that
\[
 \sum_{x\in\Z} p_x \left (\theta \sigma_{y,z}(\eta^a_{x-1,x})+(1-\theta) \sigma_{y,z}(\eta^c_{x-1,x})- \sigma_{y,z}(\eta)  \right)
=  p_y D^-_y \sigma_{y,z}(\eta) + p_z D^-_z \sigma_{y,z}(\eta)
\]
and hence
\[
\mathcal{L} \sigma_{y,z}(\eta) = \left(L_y+L_z\right) \sigma_{y,z}(\eta).
\]
Using this in (\ref{eq:generator_separate}) completes the proof. \qed
\\
\\
\begin{Corollary}
For $y=(y_1,\ldots,y_{2n})$ with $y_1< \ldots < y_{2n}$, the expectation 
$ u(t,y,\eta) = \E_{\eta} \left[ \Sigma_{y}(\Eta_t) \right]$ satisfies a system of linear
differential equations,
\begin{eqnarray}\label{corollary:eqn}
\partial_t u(t,y,\eta) =\sum_{i=1}^{2n}\left(q_x D^{+}_{y_i} + p_{x} D^{-}_{y_i}\right) u(t,y,\eta)
\end{eqnarray}
\end{Corollary}
\textbf{Proof.}
For $y=(y_1,\ldots,y_{2n})$ with $y_1< \ldots < y_{2n}$, the expectation 
$ u(t,y,\eta) = \E_{\eta} \left[ \Sigma_{y}(\Eta_t) \right]$ satisfies
\[
\partial_t u(t,y,\eta) =   \mathcal{L} u(t,y,\eta) = \E_{\eta} \left[ \mathcal{L} \Sigma_{y}(\Eta_t) \right]
= \E_{\eta} \left[ \sum_{i=1}^{2n} \, L_{y_i} \Sigma_{y}(\Eta_t) \right]
= \sum_{i=1}^{2n} \, L_{y_i} u(t,y,\eta).
\]
The first equality is the Kolmogorov equation for the Markov process $\Eta$, the second equality is due to the
Markov semigroup and its generator commuting, the third equality is due to Lemma \ref{lemma:generator}. The statement of the corollary follows using the explicit expression
(\ref{eq:particlegenerator}) for the operator $L$.
\qed
}
\\
\\
Thus the function $u(t,y,\eta)$ allows us to recast an infinite dimensional Kolmogorov equation in $(t,\eta)$ as a finite dimensional ODE in $(t,y)$.  
The next lemma shows that this ODE is exactly solved by a Pfaffian built
out of $u(t, (y_1, y_2), \eta)$'s.
The uniqueness of solutions to the system of ODE's (\ref{corollary:eqn})
then will allow us to derive a Pfaffian expression
for the expectation of the  duality functions.
\begin{Lemma} \label{lemma:scalarpfaffian}
For all $\eta \in \{0,1\}^{\Z}$,  
for all $n \geq 1$, $y_1 \leq \ldots \leq y_{2n}$ and $t \geq0$, 
\[
\E_{\eta} \left[ \Sigma_{y}(\Eta_t) \right] = \Pf(K^{(2n)}(t,y))
\]
where $K^{(2n)}(t,y)$
is the anti-symmetric $2n \times 2n$ matrix with entries $K_t(y_i,y_j)$ for $i <j$ for the function
$K_t$ defined in (\ref{eq:scalarkernel}), that is $K_t(y,z) = \E_{\eta} [ \sigma_{y,z}(\Eta_t) ]$.
\end{Lemma}
\textbf{Proof.} 
For $n \geq 1$ denote
\begin{eqnarray*}
V_{2n} &=& \{y \in \Z^{2k}: y_1< \ldots < y_{2n}\}, \\
\overline{V}_{2n} &=& \{y \in \Z^{2k}: y_1\leq \ldots \leq y_{2n}\}, \\
\partial V_{2n}^{(i)} & = & \{y \in \Z^{2n}: y_1<\ldots<y_i = y_{i+1} < y_{i+2} < \ldots y_{2n} \}, \\
\partial V_{2n} & =  & \cup_{i=1}^{2n-1} \partial V_{2n}^{(i)}.
\end{eqnarray*}
We now detail a system of ODEs indexed by  $V_{2n}$, which will involve driving terms
indexed by $\partial V_{2n}$.
Fix an initial condition $\eta \in \{0,1\}^{\Z}$ and $n \geq 1$ and define 
\[
u^{(2n)}(t,y) = \E_{\eta} \left[ \Sigma_{y}(\Eta_t) \right] \qquad \mbox{for $t \geq 0
$, and  $y \in \overline{V}_{2n}$.}
\]
For $n \geq 1$, $u^{(2n)}$ solves the following system of ODEs:
\[
(\mbox{ODE})_{2n} \hspace{.3in} \left\{ \begin{array}{rcll}
\partial_t u^{(2n)}(t,y) &=& \sum_{i=1}^{2n} \, L_{y_i} u^{(2n)}(t,y) & \mbox{on $[0,\infty) \times V_{2n}$,} \\
u^{(2n)}(t,y) & = & u^{(2n-2)}(t,y^{i,i+1}) &  \mbox{on $[0,\infty) \times \partial V_{2n}^{(i)}$,} \\
u^{(2n)}(0,y) & = & \Sigma_y(\eta) &  \mbox{on $V_{2n}$.} 
\end{array} \right.
\] 
The notation $y^{i,i+1}$ is for the vector $y$ with coordinates $y_i$ and $y_{i+1}$ removed. Thus, when $n \geq 2$,
for $y \in \partial V_{2n}^{(i)}$ we have $y^{i,i+1} \in V_{2n-2}$.  
$\mbox{(ODE)}_{2n}$ is a system of ODEs indexed over $V_{2n}$. For $n \geq 2$, to evaluate $L_{y_i} u^{(2n)}$ one may need the values
of $u^{(2n)}$ at some points $y \in \partial V_{2n}^{(i)}$, which then act as driving functions for the 
differential equation.  The second equation, which we call the boundary condition, 
states that these can be deduced  from the values of $u^{(2n-2)}$. Indeed the boundary condition follows simply from the fact that
\[
\Sigma_y(\eta) = \Sigma_{y^{i,i+1}}(\eta) \qquad \mbox{for $y \in \partial V_{2n}^{(i)}$ and $n \geq 2$.}
\] 
By setting $u^{(0)}=1$, we may suppose the equation holds also for $n=1$, encoding the fact that $u^{(2)}(t,(y,y))=1$
for all $t \geq 0, y \in \Z$. 

The infinite sequence of equations $(\mbox{(ODE)}_{2n}: n=1,2,\ldots)$ are uniquely solvable, within the class of 
continuously differentiable functions satisfying $\sup_{t \geq 0} \sup_{y \in V_{2n}} |u^{(2n)}(t,y)| < \infty$.
Indeed the boundary condition for $u^{(2)}$ is simply that $u^{(2)}(t,(y,y)) = 1$, and standard (weighted) 
Gronwall estimates show uniqueness of solutions of $\mbox{(ODE)}_{2}$. Inductively, the boundary condition for $u^{(2n)}$ is given 
by the uniquely determined values of $u^{(2n-2)}$ and hence $u^{(2n)}$ can be found uniquely from $\mbox{(ODE)}_{2n}$.

We now check that $(\Pf(K^{(2n)}(t,y)) :n=1,2,\ldots)$ also satisfies $(\mbox{(ODE)}_{2n}: n=1,2,\ldots)$.   
First we consider the initial conditions. Fix $y \in V_{2n}$ and choose $y_0 \leq y_1$. For $\theta>0$  the entries
in the Pfaffian at time zero can be rewritten as
\[
K_0(y_i,y_j) = \sigma_{y_i,y_j}(\eta) = \frac{(-\theta)^{\eta[y_0,y_j)}}{(-\theta)^{\eta[y_0,y_i)}}.  
\]
The Pfaffian identity (see appendix)
\begin{equation} \label{eq:quotientpfaffian}
\Pf \left(\frac{a_i}{a_j}:1 \leq i < j \leq 2n \right) = \frac{a_1 a_3, \ldots a_{2n-1}}{a_2 a_4 \ldots a_{2n}} \quad \mbox{when $a_i \neq  0$ for all $i$,}
\end{equation}
shows that
\[
\Pf (K^{(2n)}(0,y)) = \prod_{i=1}^n \frac{(-\theta)^{\eta[y_0,y_{2i-1})}}{(-\theta)^{\eta[y_0,y_{2i})}} = \Sigma_y(\eta).
\]
By letting $\theta \downarrow 0$ we find the same final identity is true when $\theta=0$.

Next we check the boundary conditions. We fix $y \in \partial V_{2n}^{(i)}$, $t \geq 0$ and write $K^{(2n)}$ for $K^{(2n)}(t,y)$. 
By conjugating with a suitable elementary matrix $E$, the matrix
\[
\hat{K}^{(2n)}= E^T K^{(2n)} E
\]
is the result of subtracting row $i+1$ from row $i$, and column $i+1$ from column $i$. The Pfaffian identity
$\Pf(E^T  A E) = \det(E) \Pf(A)$ ensures that $\Pf(\hat{K}^{(2n)}) = \Pf(K^{(2n)})$. However the equality
$y_i=y_{i+1}$ implies that the $i$th row of $\hat{K}^{(2n)}$ has all zero entries except for $\hat{K}^{(2n)}_{i \,i+1}= 1$. Performing a Laplace expansion
(see appendix) of the Pfaffian of $\hat{K}^{(2n)}$ along row $i$ shows that, when $n \geq 2$, 
\[
\Pf(\hat{K}^{(2n)}(t,y)) = \Pf(K^{(2n-2)}(t,z)) \quad \mbox{where $z=y^{i,i+1}$.}
\] 
When $n=1$ we obtain that $\Pf(\hat{K}^{(2)})=1$. This is exactly the desired boundary condition.

Finally we check the differential equation in $\mbox{(ODE)}_{2n}$. The entries in $K^{(2n)}(t,y)$ solve $(\mbox{ODE})_2$, 
that is $\partial_t K_t(y_i,y_j) = (L_{y_i} + L_{y_j}) K_t(y_i,y_j)$. 
The Pfaffian $\Pf(K^{(2n)}(t,y))$ is a sum of terms each of product form
\begin{equation} \label{eq:pfaffianterm}
K_t(y_{\pi_1},y_{\pi_2}) \ldots K_t(y_{\pi_{2n-1}},y_{\pi_{2n}})
\end{equation}
for some permutation $\pi$, containing each of the variables $(y_i:i \leq 2n)$ exactly once. 
Hence each term, and therefore the entire Pfaffian, 
solves the desired equation $\partial_t u = \sum_{i=1}^{2n} \, L_{y_i} u$ when $y \in V_{2n}$.

Note that $|K_t(x,y)| \leq 1$ and hence the Pfaffian $\Pf(K^{(2n)}(t,y))$ is a uniformly bounded function on
$[0,\infty) \times V_{2n}$. Uniqueness of solutions to the sequence ($(\mbox{(ODE)}_{2n}: n=1,2,\ldots)$ now implies that
\[
\E_{\eta} \left[ \Sigma_{y}(\Eta_t) \right] = u^{(2n)}(t,y) = \Pf(K^{(2n)}(t,y)), \quad \mbox{for $n \geq 1, t \geq 0$ and $y \in V_{2n}$.} 
\]
This lemma states that this identity holds also for $y \in \overline{V}_{2n}$. However, by repeating the argument for 
the boundary conditions, for $y \in \overline{V}_{2n}$ any equalities in $y_1 \leq y_2 \leq \ldots \leq y_{2n}$ can be removed, pair by pair, until
\[
u^{(2n)}(t,y) = u^{(2m)}(t,z), \quad \mbox{and} \quad \Pf(K^{(2n)}(t,y)) = \Pf(K^{(2m)}(t,z))
\] 
for some $m \leq n$ and $z \in V_{2m}$ , and hence equality also holds on the larger set. 
\qed
\hspace{.1in}

\textbf{Proof of Theorem \ref{T1}.} The correlation functions $\E[\Eta_t(x_1) \ldots \Eta_t(x_n)]$ can be recovered from the product spin expectations.
Indeed
\begin{eqnarray*}
D^+_z \sigma_{y,z}(\eta)  &=&  \sigma_{y,z+1}(\eta) - \sigma_{y,z}(\eta) \\
&=&  \sigma_{y,z}(\eta) \left( (-\theta)^{\eta(z)} -1 \right) \\
& = & - (1+\theta) \, \eta(z) \, \sigma_{y,z}(\eta)
\end{eqnarray*}
so that
\begin{equation} \label{eq:diffsigma}
\left. D^+_z \sigma_{y,z}(\eta) \right|_{z=y}  = - (1+\theta) \, \eta(y).
\end{equation}
Thus out of a spin pair, we can reconstruct a single occupancy variable by first a discrete derivative, and then an evaluation. 
We now iterate this to get at multiple occupancy variables. Fix $n \geq 1$ and
consider $y=(x_1,\hat{x}_1, \ldots, x_n, \hat{x}_n) \in \overline{V}_{2n}$ where
 \[
x_1 \leq \hat{x}_1 < x_2 \leq \hat{x}_2 < x_3 \leq \ldots<x_n \leq \hat{x}_n.
\]
The restriction that $\hat{x}_i < x_{i+1}$ allows us to 
apply the operators $D^+_{\hat{x_1}}, \ldots, D^+_{\hat{x_n}}$ to both sides of the   
identity (from Lemma  \ref{lemma:scalarpfaffian})
\[
\E_{\eta} \left[ \Sigma_{y}(\Eta_t) \right] = \Pf(K^{(2n)}(t,y)).
\]
The left hand side becomes
\[
(-1)^n (1+\theta)^n \E_{\eta}  \left[ \Eta_t(x_1) \ldots \Eta_t(x_n) \prod_{i=1}^n \sigma_{x_i,\hat{x}_i}(\Eta_t) \right].
\]
After setting $x_i = \hat{x_i}$ for all $i$ we reach the correlation $(-1)^n (1+\theta)^n \E_{\eta} [ \Eta_t(x_1) \ldots \Eta_t(x_n)]$.
Applying the operators $D^+_{\hat{x_1}}, \ldots, D^+_{\hat{x_n}}$ to the Pfaffian on the right hand side
preserves the Pfaffian structure. Indeed applying $D^+_{\hat{x}_1}$ to the single Pfaffian term (\ref{eq:pfaffianterm})
will change only a single factor in the product, namely
\[
K_t(y_{\pi_{2i-1}},y_{\pi_{2i}}) \to \left\{ \begin{array}{ll}
D^+_1 K_t(y_{\pi_{2i-1}},y_{\pi_{2i}}) & \mbox{if $\pi_{2i-1}=1$,} \\
D^+_2 K_t(y_{\pi_{2i-1}},y_{\pi_{2i}}) & \mbox{if $\pi_{2i}=1$.}
\end{array} \right.
\]
The terms can then be summed into a new Pfaffian where the entries in the second row and column, 
which are the only entries containing the variable $\hat{x}_1$, are changed. Repeating this for the 
operators $D^+_{\hat{x_2}}, \ldots, D^+_{\hat{x_n}}$, and then setting $x_i = \hat{x_i}$ for all $i$ we still
have a Pfaffian, and the final entries are in the form
\[
\frac{-1}{1+\theta} \left( \begin{array}{cc} 
K_t(y,z) & D^+_2K_t(y,z) \\
D^+_1K_t(y,z) & D^+_1 D^+_2 K_t(y,z)
\end{array} \right)
\]
A final conjugation with a determinant one diagonal matrix, with entries that are alternating $\pm i$,
will adjust the minus signs to give exactly by the kernel $\K$ stated in the theorem. \qed.

\vspace{.1in}
\noindent

\textbf{Remarks.}

\noindent
\textbf{1.} For random initial conditions the point process $\eta_t$ will not in general be a Pfaffian point process, although by 
conditioning on the initial condition, the intensities can always be written as the expectation of a Pfaffian. However,
under some random initial conditions, $\eta_t$ does remain a Pfaffian point process. 
Indeed, examining the proof of Lemma \ref{lemma:scalarpfaffian} , one needs only that the expectation
$\E \left[ \Sigma_y(\Eta_0) \right]$, for $y \in V_{2n}$,
can be written as a Pfaffian $\Pf(\Phi(y_i,y_j): i < j)$ for some $\Phi:V_2 \to \R$. 
One then replaces the initial condition in the equation (\ref{eq:ode2}) for $K_t(y,z)$ by $\Phi(y,z)$
and the rest of the argument goes through. 
An example is where the sites $(\Eta_0(x): x \in \Z)$ are independent with $\eta_0(x)$ a Bernoulli\,($B_x$) variable, 
then 
\[
\E \left[ \Sigma_y(\Eta_0) \right] 
= \Pf \left( \prod_{k \in [y_i,y_j)}(1-(1+\theta)B_k): i < j\right).
\]

\noindent
\textbf{2.} A slightly more combinatorial way of writing out the argument for the last part of the proof of Theorem \ref{T1}
is as follows. 
Starting from (\ref{eq:reconstruct}) we may reconstruct the product intensities as 
\begin{eqnarray} \label{eq:correlation}
 \E \left[ \prod_{i=1}^n \Eta_t(x_i) \right] 
 &=& \E \left[ \prod_{i=1}^n \frac{(1- \sigma_{x_i,x_i+1}(\Eta_t))}{1+\theta} \right] \nonumber \\
 &=& (1+\theta)^{-n} \sum_{m=0}^n(-1)^m \sum_{\stackrel{y_1<\dots<y_m}{\in\{x_1,\dots,x_n\}}} 
\E \left[\prod_{i=1}^m \sigma_{y_i,y_i+1}(\Eta_t) \right].
\end{eqnarray}
Since the vector $y^{(2m)}=(y_1,y_1+1,\ldots,y_m,y_m+1) \in \overline{V}_{2m}$ we may apply Lemma \ref{lemma:scalarpfaffian} to see that
\[
 \E \left[ \prod_{i=1}^n \Eta_t(x_i) \right] 
=  (1+\theta)^{-n} \sum_{m=0}^n(-1)^m \sum_{\stackrel{y_1<\dots<y_m}{\in\{x_1,\dots,x_n\}}} \Pf \left(K^{(2m)}(t,y^{(2m)})\right).
\]
This sum may be recombined as the single Pfaffian
\[
 \E \left[ \prod_{i=1}^n \eta_t(x_i) \right] 
=  (-1)^n (1+\theta)^{-n} \Pf \left(K^{(2n)}(t,y^{(2n)}) - J_{2n} \right)
\]
where $J_{2n}$ is the block diagonal matrix formed by $n$ copies of {\tiny $\left(\! \! \! \!\begin{array}{rl} 0 & 1\\ -1&0 \end{array} \! \! \right)$}.
Indeed  this is a special case of the general formula for $\Pf(A+B)$ (see (\ref{eq:Jcase}) in the appendix). This shows 
that $\Eta_t$ is a Pfaffian point process
with the kernel 
$\tilde{\K}$ given, 
for $y<z$, by
\[
\tilde{\K}(y,z) = \frac{-1}{1+\theta} \left( \begin{array}{cc} 
K_t(y,z) & K_t(y,z+1) \\
K_t(y+1,z) & K_t(y+1,z+1)
\end{array} \right),
\]
and $\tilde{\K}_{12}(y,y) =  \frac{-1}{1+\theta} (K_t(y,y+1)-1)$, and other entries determined by the symmetry conditions.
That both kernels $\K$ and $\tilde{\K}$ determine the same point process can be seen by a simple transformation: conjugation by an elementary matrix which subtracts the first row and column from the second row and column,
followed by a conjugation to adjust the positions of the minus signs as before, will transform the entries of
$\tilde{\K}$ into those for $\K$.
\section{Examples} \label{s3}
We prove the various continuum limits given as the four examples (A),(B),(C) and (D) listed in Theorem \ref{T2}. The fact that the discrete kernels satisfy two dimensional discrete heat equations make all the examples easy to guess, and they are not that 
much harder to prove. We divide our proof into three steps: (1) convergence of the lattice particle systems to their
continuum analogues for finite particle systems; (2) establishing the Pfaffian kernel for these finite continuum system
by limits of the lattice Pfaffian kernels; (3) taking limits to obtain the continuum systems under maximal entrance 
laws and solving explicitly for the kernel. 
(In \cite{TZ:11} we argued directly with the continuum Brownian models, and we don't doubt that this would also be possible 
for these four examples.)

We use the space $\mathcal{M}_{LFP}$ of locally finite point measures on $\R$ with the topology
of vague convergence, which has simply checked compactness properties. 
Our point processes can then be thought of taking values in the measurable subset 
$\mathcal{M}_0 \subset \mathcal{M}_{LFP}$ of simple measures, that is where all atoms have mass one. 

\textbf{Step 1. Finite particle continuum approximation.}
Working first with finite particle systems avoids complicated weak convergence arguments for
infinite systems. We fix $N \geq 1$ and  initial particle positions $a_1<a_2< \ldots < a_N$. 
In example (C) we assume $a_1>0$ and in example $(D)$ we assume $a_1 \geq 0$.
We choose finite lattice initial conditions
\begin{equation} \label{ics}
\eta_0(x) = \left\{ \begin{array}{ll}
1 & \mbox{if $x = a^{(\epsilon)}_i$ for $i = 1,\ldots,N$,} \\
0 & \mbox{otherwise,}
\end{array} \right.
\end{equation}
where $\epsilon a^{(\epsilon)}_i \to a_i$ as $\epsilon \to 0$ for $i =1,\ldots,N$. 
We fix $t>0$ and set 
$ X^{(\epsilon)}_t(dx) = \sum_{x \in \Z}  \eta_{\epsilon^{-2}t}(x) \delta_{\epsilon x}.$
The conclusion we want is that
 \[
 (X^{(\epsilon)}_t: t \in [0,T]) \to  (X_t: t \in [0,T]) 
 \]
in the Skorokhod space $D([0,T], \mathcal{M}_{LFP})$. One way to check this
is via the continuous mapping principle (see \cite{billingsley}), from the weak convergence of the scaled non-interacting 
simple random walk paths $(Z^{{\epsilon},1},\ldots,Z^{{\epsilon},N}) \to (B^{(1)}, \ldots, B^{(N)})$, in
$D([0,T],\R^N)$, to independent Brownian motions. For the purely annihilating system
there is a deterministic map $F: D([0,T],\R^N) \to D([0,T],\mathcal{M}_{LFP})$ that 
'prunes' the paths appropriately at collision times, satisfying 
\[
(X^{(\epsilon)}_t: t \in [0,T]) = F(Z^{{\epsilon},1}_t,\ldots,Z^{{\epsilon},N}_t:t \in [0,T]).
\]
We will not detail the natural definition of $F$, but it is straightforward when the collision times
of the $N$ paths are all distinct, and and when two collision times are equal we may define $F$ to be 
the zero function. The law of the limiting Brownian paths does not charge the 
discontinuities of $F$, and desired result follows. 

For mixed models the argument can be adjusted to allow for the random reaction that is required, 
and we give an informal description.
From $N$ initial particles there are at most $N$ reactions. We can fix a vector of reactions
$r=(r_1,\ldots,r_N)$ where $r_i \in \{A,C^+,C^-\}$ reflect the successive 
decisions: annihilate, coalesce onto the 'higher' of the two colliding paths, coalesce onto the 'lower' of the two colliding paths.
For a fixed $r$ there is again a deterministic map $F_r: D([0,T],\R^N) \to D([0,T],\mathcal{M}_{LFP})$ 
that prunes the $N$ paths according to the set of decisions (not all the decisions may be needed by time $T$).
Then, for a bounded and continuous $\Phi: D([0,T],\mathcal{M}_{LFP}) \to \R$
\[
 E[ \Phi(X^{(\epsilon)}_t: t \in [0,T])] = \sum_r \theta^{A(r)} (1-\theta)^{N-A(r)} 
 E[\Phi(F_r(Z^{{\epsilon},1}_t,\ldots,Z^{{\epsilon},N}_t:t \in [0,T]))]
\]
where $A(r)$ is the number of annihilations in $r$, and the result follows as for the coalescing case.
For example $(D)$ we apply the same to reflected paths, and 
for example $(C)$ we include the times when a particle first hits zero as an extra collision, which 
reduces the number of particles by one but which requires no random reaction.

\textbf{Step 2. From lattice kernels to continuum kernels}
We record the limiting continuum Pfaffian kernels, for systems with  finite initial conditions.
\begin{Proposition} \label{P1}
The continuum models with finite initial conditions $\eta_0 = \sum_{i=1}^{N} \delta_{a_i} $ are 
Pfaffian point processes at a fixed $t>0$ with kernel in the form (\ref{ctmkernel})
where $K_t^{(c)}(y,z)$ is the unique solution to the following PDES:
\begin{itemize}
\item \textbf{(A),(B) Coalescing/annihilating Brownian motions.}
\begin{equation}
\label{Bpde}
\left\{  \begin{array}{rcll} 
\partial_t K^{(c)}_t(y,z) &=& \Delta K^{(c)}_t(y,z) & \mbox{for $y<z, \, t>0$,}\\
K^{(c)}_t(y,y) & = & 1 & \mbox{for $y \in \R, \, t>0$,} \\
K^{(c)}_0(y,z) & = & (-\theta)^{\eta_0[y,z)} , &  \mbox{for $y \leq z$,}
\end{array} \right.
\end{equation}
\item \textbf{(C) Coalescing/annihilating Brownian motions on $[0,\infty)$ killed at $\{0\}$.}
\begin{equation}
\label{BKpde}
\left\{
\begin{array}{rcll} 
\partial_t K^{(c)}_t(y,z) &=& \Delta K^{(c)}_t(y,z) & \mbox{for $0<y<z, \, t>0$,}\\
K^{(c)}_t(y,y) & = & 1 & \mbox{for $0<y, \, t>0$,} \\
\partial_1 K^{(c)}_t(0,z) &=& 0, & \mbox{for $0<z, \, t>0$,} \\
K^{(c)}_0(y,z) & = & (-\theta)^{\eta_0[y,z)}  &  \mbox{for $0 \leq y \leq z$.}
\end{array} \right.
\end{equation}
\item \textbf{(D) Coalescing/annihilating reflected Brownian motions on $[0,\infty)$.}
\begin{equation}
\label{BRpde}
\left\{ 
\begin{array}{rcll} 
\partial_t K^{(c)}_t(y,z) &=& \Delta K^{(c)}_t(y,z) & \mbox{for $0<y<z, \, t >0$,}\\
K^{(c)}_t(y,y) & = & 1 & \mbox{for $y >0, \, t>0$,} \\
K^{(c)}_t(0,z) &=&  \Phi_t(z) & \mbox{for $0<z, t >0$,} \\
K^{(c)}_0(y,z) & = & (-\theta)^{\eta_0[y,z)}, & \mbox{for $0 \leq y \leq z$,}
\end{array} \right.
\end{equation}
where $(\Phi_t(z): z ,t \geq 0)$ solves
\[
\left\{ 
\begin{array}{rcll} 
\partial_t \Phi_t(z) &=& \Delta \Phi_t(z) & \mbox{for $ 0<z, \, t>0$,}\\
\Phi_t(0) & = & 1, & \mbox{for $t>0$.} 
\end{array} \right.
\]
\end{itemize}
\end{Proposition}
In each case we will show the convergence of the approximating lattice system $X^{(\epsilon)} \to X$. 
The following lemma 
(whose proof is at the end of the appendix) is a natural approximation lemma 
for lattice kernels to continuum kernels.  
\begin{Lemma} \label{kernelconvergence}
For $\epsilon>0$, let  $X^{(\epsilon)}(dx)$ be random point measure on $\R$ whose atoms
from a  Pfaffian point process on $\epsilon\Z$ with kernel $\K^{(\epsilon)}$. 
Suppose that
\begin{equation}
\sup_{\epsilon>0} \| \epsilon^{-1} \K^{(\epsilon)}\|_{\infty} = 
\sup_{\epsilon>0}
\sup_{y,z,i,j} \epsilon^{-1} \, |\K^{(\epsilon)}_{ij}(y,z)|  < \infty  \label{eq:kernelbound1}
\end{equation}
and
\begin{eqnarray}
&& \lim_{\epsilon\downarrow0} \epsilon^{-1} \K_{ij}^{(\epsilon)} (y_\epsilon,z_\epsilon) = \K^{(c)}_{ij} (y,z), 
\quad \mbox{for $i,j \in\{1,2\},$}\label{eq:kernelconvergence1} \\
&& \qquad \mbox{when $(y_\epsilon,z_{\epsilon}) \to (y,z)$ with $y<z$, or when $y_\epsilon = z_\epsilon \to z=y$,} \nonumber 
\end{eqnarray}
for some continuum kernel $\K^{(c)}: \R^2 \to \R^{2 \times 2}$.
Then $X^{(\epsilon)} \to X$ in distribution as $\epsilon \downarrow 0$, on the space $\mathcal{M}_{LFP}$  
and the limit $X$ is simple, and is a Pfaffian point process with kernel $\K^{(c)}$.
\end{Lemma}
Note that in our examples the limiting kernel $\K^{(c)}(x,y)$ will be discontinuous at $x=y$.

\textbf{(A),(B)}.  This is the case $p_x=q_x=1$ for all $x \in \Z$.
The corresponding one-particle generator is then the discrete Laplacian $\Delta = D^+ + D^-$.
Lemma \ref{lemma:generator} shows that 
the entries of the Pfaffian kernel $\K^{(\epsilon)}$ for $X^{(\epsilon)}$ are given for $y<z$ in $\epsilon\Z$ by 
\[
\begin{array}{ccl}
 \K^{(\epsilon)}(y,z) 
 &=& \frac{1}{1+\theta} \left( \begin{array}{cc} 
K_{\epsilon^{-2}t}(\epsilon^{-1}y,\epsilon^{-1}z) & - D^+_2K_{\epsilon^{-2}t}(\epsilon^{-1}y,\epsilon^{-1}z) \\
- D^+_1K_{\epsilon^{-2}t}(\epsilon^{-1}y,\epsilon^{-1}z) & D^+_1 D^+_2 K_{\epsilon^{-2}t}(\epsilon^{-1}y,\epsilon^{-1}z)
\end{array} \right),
\end{array}
\]
and by $\K^{(\epsilon)}_{12}(y,y) =  \frac{-1}{1+\theta} \; D^+_2 K_{\epsilon^{-2}t}(\epsilon^{-1}y,\epsilon^{-1}y)$.
We rescale the scalar kernel by defining
\[
K^{(\epsilon)}_t(y,z) = K_{\epsilon^{-2}t}(\epsilon^{-1}y,\epsilon^{-1}z) \quad \mbox{for $y,z \in \epsilon \Z$.}
\]
For  $f:\epsilon\Z \to \R$ we set
\[
D^{(\epsilon),+} f(x) = \frac{f(x+\epsilon)-f(x)}{\epsilon}, \quad
 \Delta^{(\epsilon)}f(x) = \frac{f(x+\epsilon) + f(x-\epsilon) - 2f(x)}{\epsilon^2}.
\]
Then
\begin{equation} \label{latticekernelform}
\begin{array}{ccl}
 \K^{(\epsilon)}(y,z) 
 &=& \frac{1}{1+\theta} \left( \begin{array}{cc} 
K^{(\epsilon)}_t(y,z)  & - \epsilon D^{(\epsilon),+}_2 K^{(\epsilon)}_t(y,z)  \\
- \epsilon D^{(\epsilon),+}_1 K^{(\epsilon)}_t(y,z)  & \epsilon^2 D^{(\epsilon),+}_1 D^{(\epsilon),+}_2  K^{(\epsilon)}_t(y,z) 
\end{array} \right),
\end{array}
\end{equation}
and $ K^{(\epsilon)}_t(y,z)$ solves, for $y,z \in \epsilon \Z$,
\begin{equation} \label{latticepde}
\left\{  \begin{array}{rcll} 
\partial_t K^{(\epsilon)}_t(y,z) &=& \left( \Delta^{(\epsilon)}_y + \Delta^{(\epsilon)}_z \right) K^{(\epsilon)}_t(y,z) & \mbox{for $y<z, t>0$,}\\
K^{(\epsilon)}_t(y,y) & = & 1, & \mbox{for all $y, t>0$,} \\
K^{(\epsilon)}_0(y,z)& = & (-\theta)^{\eta_0[\epsilon^{-1}y, \epsilon^{-1}z)} & \mbox{for $y \leq z$.}
\end{array} \right.
\end{equation}
By conjugation with a diagonal matrix with diagonal entries $\epsilon^{1/2}$ and $\epsilon^{-1/2}$ an alternative kernel for 
$X^{(\epsilon)}$, in the right form for Lemma \ref{kernelconvergence}, is
\[
 \tilde{\K}^{(\epsilon)} (y,z) = \frac{\epsilon}{1+\theta} \left( \begin{array}{cc} 
K^{(\epsilon)}_t(y,z)  & - D^{(\epsilon),+}_2 K^{(\epsilon)}_t(y,z)  \\
- D^{(\epsilon),+}_1 K^{(\epsilon)}_t(y,z)  & D^{(\epsilon),+}_1 D^{(\epsilon),+}_2  K^{(\epsilon)}_t(y,z) 
\end{array} \right).
\]
Checking the hypotheses (\ref{eq:kernelbound1}) and (\ref{eq:kernelconvergence1}) amounts to checking that the lattice 
approximations to the two dimensional continuum PDE (\ref{Bpde}) converge uniformly, at a fixed $t>0$, 
along with their first and second derivatives. The required estimates are quite standard and 
we omit the proof here, and also for the examples C,D below. Some details (however for different initial conditions) 
are contained in the thesis \cite{thesis}.

\textbf{(C).} We take
\begin{equation} \label{killedparameters}
q_x = \left\{ \begin{array}{ll}
1 & \mbox{for $x \geq 1$,} \\ 
2 & \mbox{for $x=0$,} \\
0 & \mbox{for $x <1$,}
\end{array} \right.
 \qquad 
p_x = \left\{ \begin{array}{ll}
1 & \mbox{for $x \geq 1$,} \\ 
0 & \mbox{for $x <1$.}
\end{array} \right.
\end{equation}

No particle ever visits $\{\ldots,-3,-2\}$, and particles that reach $\{-1\}$ never escape. So we
restrict attention to the point process $(\Eta_t(x): x=0,1,\ldots)$ and define
$X^{(\epsilon)}(dx) = \sum_{x \geq 0}  \eta_{\epsilon^{-2}t}(x) \delta_{\epsilon x}$ as a measure on
$\epsilon \Z_+ = \{ 0, \epsilon, 2 \epsilon, \ldots\}$. A single particle acts as a simple random walk on
$\{0,1,\ldots\}$ with a certain rate of being killed whenever it is at zero. Under diffusive rescaling this
process becomes a Brownian motion that is instantly killed at the origin (which allows the treatment 
for interacting particles to go through as in step 1.) 

The reason for choosing $q_0=2$ is so that the corresponding one particle generator $L$ defined in (\ref{eq:particlegenerator}) 
is the generator for a reflected 
random walk on $\{0,1,\ldots\}$, which
jumps $x \to x \pm 1$ at rate $1$ for $x \geq 1$, and also jumps $0 \to 1$ at rate $2$.
This can be realised as the absolute value of a simple random walk on $\Z$, and this helped us in some 
of the estimates showing the lattice PDE converged to the continuum PDE. 

Then $X^{(\epsilon)}$ is Pfaffian with kernel of the form
(\ref{latticekernelform}) where Lemma \ref{lemma:generator} shows that $ K^{(\epsilon)}_t(y,z)$ solves (\ref{latticepde}) for 
$\epsilon \leq y < z$. 
As expected, the Neumann boundary condition emerges in the limiting 
continuum PDE (\ref{BKpde}) when $\epsilon \to 0$.
 
\textbf{(D).} 
To obtain reflected random walks on $\Z_+$ we take
\[
q_x = \left\{ \begin{array}{ll}
1 & \mbox{for $x \geq 1$,} \\ 
0 & \mbox{for $x <1$,}
\end{array} \right.
 \qquad 
p_x = \left\{ \begin{array}{ll}
1 & \mbox{for $x \geq 1$,} \\ 
0 & \mbox{for $x <1$.}
\end{array} \right.
\]
We may restrict attention to the process $(\Eta_t(x): x \geq 0)$. Note that the 
corresponding one particle generator $L$ defined in (\ref{eq:particlegenerator}) is the generator 
for simple random walk on $\{0,1,\ldots\}$ absorbed at $0$. 
Thus $X^{(\epsilon)}$ is Pfaffian with kernel of the form
(\ref{latticekernelform}) where $ K^{(\epsilon)}_t(y,z)$ solves (\ref{latticepde}) for 
$\epsilon \leq y < z$. The boundary condition can be found by examining the 
expectation $K_t(0,z) = \E_{\eta} \left[ \sigma_{0,z}(\Eta_t) \right]$. Applying the generator
we find that
\[
\partial_t K_t(0,z) = L_z K_t(0,z) \quad \mbox{for $z \geq 1$,}
\]
with the boundary condition $K_t(0,0)=1$. Under scaling the limiting boundary condition 
becomes $\Phi_t(z)$ as stated in example D of Theorem \ref{T2}. 

\textbf{Step 3. Maximal entrance laws.}
Infinite systems of coalescing particles are easy to build by adding one particle at a time. 
Infinite systems of annihilating particles perhaps require a bit more care. In \cite{TZ:11}
a pathwise construction is bypassed since only fixed time properties are studied, and instead
this paper uses a Feller transition kernel $p_t(\mu,d \nu)$ on the measurable subset $\mathcal{M}_0$
of simple measures within the space $\mathcal{M}_{LFP}$. Our mixed models also have Feller 
transition kernels, and we can construct these by exploiting the Pffafian structure.
Suppose $\eta^{(N)} \to \eta \in \mathcal{M}_0$, where $\eta^{(N)}$ have finitely many particles. 
Then 
\[
(-\theta)^{\eta^{(N)}[y,z)} \to (-\theta)^{\eta[y,z)} \quad \mbox{for almost all $(y,z)$.}
\]
Indeed this holds for all $y,z$ that are not the position of atoms in $\eta$. These
functions are the initial conditions for the pdes (\ref{Bpde}), (\ref{BKpde}), (\ref{BRpde}) 
that determine the Pfaffian kernels. The solutions to these pdes then converge, at a fixed $t>0$, together
with their first and second derivatives. Thus the associated Pfaffian kernels converge
(in a bounded pointwise manner) and Lemma \ref{kernelconv2} in the appendix 
states this is sufficient for the associated point processes $\eta^{(N)}_t$ to converge in law to a limit.
This defines a kernel $p_t(\eta,d \nu)$ for all $\eta \in \mathcal{M}_0$. 
To check the semigroup property it is sufficient to have the Feller property, and pass to the limit from the 
semigroup property for finite systems. The Feller property follows however from the same 
argument: the convergence $\eta^{(N)} \to \eta \in \mathcal{M}_0$ implies, via convergence of the 
Pfaffian kernels, the weak converge
of the laws $p_t(\eta^{(N)},d \nu) \to p_t(\eta,d \nu)$. 

To construct the maximal entrance laws we pick any sequence $\eta^{(N)}$ so that
\begin{equation} \label{conics}
(-\theta)^{\eta^{(N)}[y,z)}  \to 0 \quad \mbox{as $N \to \infty$.}
\end{equation}
When $\theta \in [0,1)$ it is easy to ask for pointwise convergence, but when $\theta=1$ this is impossible. Instead
we ask that the convergence holds in distribution (either on $\R$ or on $[0,\infty)$ for the 
killed and reflected models). If for example the atoms of $\eta^{(N)}$ are
at $k/N$, where $k=-N^2, -N^2+1,\ldots N^2-1,N^2$, then for any integrable
function $f$ on $\R$,
\[
\lim_{N\rightarrow \infty} \int_a^b dz (-1)^{\eta^{(N)}[a,z)} f(z)=0,
\]
which can be verified using a suitable modification of the Riemann-Lebesgue lemma.

Convergence in distribution is still sufficient to imply that the associated PDEs, and its derivatives, converges pointwise at
a fixed time $t>0$, as can be verified by writing out the solution in terms of the initial conditions and 
the associated Green's functions. The limiting laws we denote as $p_t(\infty,d \nu)$, and again the Feller property
allows one to check that they act as entrance laws for the Markov family. 

The name maximal is natural for coalescing 
systems, while for purely annihilating systems many increasingly dense initial conditions will not converge
to $p_t(\infty, d\nu)$ - a sequence of closely positioned pairs will annihilate quickly and leave empty regions.
The Pfaffian structure makes it clear that the convergence in distribution in (\ref{conics}) is exactly the right condition to describe the domain of attraction for the maximal entrance law. The convergence 
(\ref{conics})  holds in particular for lattice initial conditions as the spacing decreases to zero, or by choosing
Poisson initial conditions of increasing intensity. 

Finally under the maximal entrance laws the associated pdes (\ref{Bpde}), (\ref{BKpde}), (\ref{BRpde}) have
zero initial conditions, and the explicit solutions listed in Theorem \ref{T2} are straightforward to derive.

\section{Appendix: Pfaffian Facts}
%
%
%
We collect here some basic facts about Pfaffians. All the following results are contained, for example, in 
Stembridge \cite{Ste:90}. 

Throughout this section $A=\{a_{ij}\}_{i,j=1}^{2n}$ is a $2n\times2n$ anti-symmetric matrix (say with complex entries).
The Pfaffian can be defined by
 \[
\Pf(A) = \frac{1}{2^n n!} \sum_{\pi \in\Pi_{2n}} \mbox{sgn}(\pi) \prod_{i=1}^n a_{\pi_{2i-1},\pi_{2i}},
 \]
 where $\Pi_{2n}$ is the set of permutations on $\{1,2,\dots,2n\}$ and $\mbox{sgn}(\pi)$ is the sign of the permutation $\pi$.

Pfaffians are related to determinants and determinant properties often have Pfaffian analogues.
Many of these results are consequences of the following conjugation formula
\[
\Pf(CAC^T)=\det(C) \, \Pf(A) \quad \mbox{for all $2n \times 2n$ matrices $C$,} 
\]
and the identity $\det(A) = \left( \Pf(A) \right)^2.$

One may decompose a Pfaffian along any row (or column) in terms of sub-matrices, analogous to the Laplace expansion of a determinant: for any $i\in \{1,\dots,2n\}$
 \[
  \Pf(A)=\sum_{j=1,j\neq i}^{2n}(-1)^{i+j+1+\chi\{j<i\}}a_{ij}\Pf(A^{(i,j)}),
 \]
where $A^{(i,j)}$ is the $(2n-2)\times(2n-2)$ sub-matrix formed by removing the $i$-th and $j$-th rows and columns from $A$. 

For two anti-symmetric matrices $A,B$, we have 
\[
 \Pf(A+B)=\sum_U(-1)^{|U|/2}(-1)^{s(U)}\Pf(A|_U)\Pf(B|_{U^c}),
\]
where the sum is over subsets $U\subseteq\{1,\dots,2n\}$ with $|U|/2\in\{1,\dots,n\}$, $s(U)=\sum_{j\in U}j$ (with $s(\emptyset)=0$), and $A|_U$ means the matrix $A$ restricted to the rows and columns of $U$.

We use the special case of this when $-B=J_{2n}$, the canonical symplectic matrix
consisting of $n$ blocks 
$\left( \!\! \! \begin{array}{cc} 0 & 1 \\ -1 & 0 \end{array} \!\! \right)$ 
down the diagonal. 
Then setting $U_x=\{2x_1-1, 2x_1, \ldots, 2 x_{m}-1, 2x_m\}$
for $x \in \{1,2,\ldots,n\}^m$ with $x_1<x_2< \ldots < x_m$, we have
$\Pf(J_{2n}|U_x) = 1$ while $\Pf(J_{2n}|U) = 0$ for all (non-empty) $U \subseteq \{1,2,\ldots,2n\}$  not of this form. 
Note that $(-1)^{s(U_x)}=(-1)^m$.  In this case, the formula reduces to
\begin{equation} \label{eq:Jcase}
 \Pf(A-J_{2n})= \sum_{m=0}^n \sum_{x} (-1)^{n-m} \Pf(A|_U),
\end{equation}
where the sum over $x$ is precisely the sum over $x$ of the above form. 

\vspace{.1in}

\textbf{Proof of the Pfaffian identity (\ref{eq:quotientpfaffian})}. 
Let $A$ be the antisymmetric matrix with entries $A_{ij} = a_i/a_j$. By conjugating 
with a suitable elementary matrix we may subtract a multiple $a_1/a_2$ of the second row and column form the first 
row and column. This produces a new matrix $\hat{A}$ with $\Pf(A) = \Pf(\hat{A})$ but where the top row
of $\hat{A}$ is now $(0, a_1/a_2,0,0,\ldots)$. By a Laplace expansion of this top row we find
\[
\Pf(A) = \Pf(\hat{A}) = \frac{a_1}{a_2} \Pf(\hat{A}^{(1,2)}) = \frac{a_1}{a_2} \Pf(A^{(1,2)}), 
\]
and  by induction over $n$ we find $\Pf(A) = (a_1a_3 \ldots a_{2n-1})/(a_2 a_4 \ldots a_{2n})$. \qed

\vspace{.1in}

\noindent
\textbf{Proof of Lemma \ref{kernelconvergence}.}
For $a<b$, the first moments
\[
\E \left[ X^{(\epsilon)}([a,b]) \right] = \sum_{x \in \epsilon\Z \cap [a,b]} K_{12}^{(\epsilon)}(x,x)
\]
are uniformly bounded by (\ref{eq:kernelbound1}), and this implies the tightness of
$(X^{(\epsilon)}: \epsilon>0)$ as elements of $\mathcal{M}_{LFP}$.
The higher factorial moments (writing $[z]_k= z(z-1) \ldots(z-k+1)$) such as
\begin{eqnarray*}
\E \left[ [X^{(\epsilon)}([a,b])]_k \right] &=& \sum_{x_1,\ldots,x_k \in  \epsilon\Z \cap [a,b]} 
\Pf( K^{(\epsilon)}(x_i,x_j):i,j \leq k) \\
& \leq & \|\epsilon^{-1} \K^{(\epsilon)}\|^{2k}_{\infty} ((b-a)+\epsilon^{-1})^k
\end{eqnarray*}
are also uniformly bounded in $\epsilon \in (0,1]$. By Fatou's lemma, the moments of any limit points satisfy
$\E \left[ [X^{(\epsilon)}([a,b])]_k \right] \leq C_k (b-a)^k$ for all $k \geq 1$. 
This implies that any limit point is a simple point process: indeed
\begin{eqnarray*}
\pr \left[ X(\{x\}) \geq 2, \, \mbox{for some $x \in [-M,M]$} \right]
& \leq & \sum_{k=-Mm}^{Mm}  \pr \left[ X([k/m,(k+1)/m)]) \geq 2 \right] \\
& \leq &  \sum_{k=-Mm}^{Mm}  \E \left[ [X([k/m,(k+1)/m)])]_2  \right] \\
& \leq &  C_2 \sum_{k=-Mm}^{Mm}  m^{-2} \to 0 \quad \mbox{as $m \to \infty$.}
\end{eqnarray*}
Then one can pass to the limit, for finite disjoint intervals $A_1, \ldots, A_m$
\begin{eqnarray*}
\E \left[ \prod_{k=1}^n X^{(\epsilon)}(A_k) \right] & = &
\sum_{x_i \in A_i \cap \epsilon\Z:\; i \leq k} \Pf( \K^{(\epsilon)}(x_i,x_j):i,j \leq k) \\ 
& \to & \int_{A_1 \times \ldots \times A_m} \Pf(\K^{(c)}(x_i,x_j):i,j \leq k) dx_1 \ldots dx_k \\
& = & \E \left[ \prod_{k=1}^n X (A_k) \right],
\end{eqnarray*}
for any limit point $X$. Indeed the convergence follows from the assumptions (\ref{eq:kernelbound1}) and
(\ref{eq:kernelconvergence1}); since $\E[X(\{x\})]=0$ for limit points the discontinuities of the
function $\mu \to \mu(A_1) \ldots \mu(A_m)$ are not charged; 
the finite higher moments give the uniform integrability that
justify the final equality.
This shows that a limit point is a Pfaffian point process with kernel $\K^{(c)}$. 
Finally it is well known (see the remark after this proof) that the fact that the kernel $\K^{(c)}$ is locally bounded is sufficient to 
determine the law of the associated
Pfaffian point process in $\mathcal{M}_{LFP}$, which implies that limit points are unique and thereby the
convergence of $X^{(\epsilon)}$. \qed

Very similar arguments establish the following kernel convergence lemma
for continuous kernels.
\begin{Lemma} \label{kernelconv2}
Suppose $(X^{(N)}(dx):N \geq 1)$ are random point measure on $\R$ whose atoms
form Pfaffian point processes on $\R$ with kernel $\K^{(N)}$. 
Suppose that
\[
\sup_{N} \|  \K^{(N)}\|_{\infty} = 
\sup_{N}
\sup_{y,z,i,j}  \, |\K^{(N)}_{ij}(y,z)|  < \infty  
\]
and 
\[
\lim_{N \to \infty}  \K_{ij}^{(N)} (y,z) = \K_{ij} (y,z), 
\quad \mbox{for $i,j \in\{1,2\},$ and $y,z \in \R$} 
\]
for some limiting kernel $\K: \R^2 \to \R^{2 \times 2}$.
Then $X^{(N)} \to X$ in distribution as $N \to \infty$, on the space $\mathcal{M}_{LFP}$  
and the limit $X$ is simple, and is a Pfaffian point process with kernel $\K$.
\end{Lemma}

\textbf{Remark.} To see that a locally bounded kernel $\K$ is sufficient to determine the law of associated point process $X$
one may argue as follows: (i) all moments of $X(f)$, where $f = \sum_{i=1}^m c_i I(A_i) \geq 0$ and $A_i$ are finite intervals, are given in terms of the kernel; 
(ii) by Hadamard's inequality
\[
|\Pf( \K(x_i,x_j):i,j \leq k)| = |\det( \K(x_i,x_j):i,j \leq k)|^{1/2} 
\leq \|\K\|_{\infty}^{k} (2k)^k; 
\]
(iii) the factorial moments $\E[[X([a,b])]_k]$ are bounded by $C(a,b,\K)^k k^k$; (iv) the moments
$\E[|X([a,b])|^k]$ are also bounded by $C'(a,b,\K)^k k^k$; (v) the moment problem for $X(f)$ is well posed. 

\noindent
{\bf Acknowledgements.} B.G. supported by EPSRC grant EP/H023364/1; M. P. and R. T. supported by 
EPSRC grant No.\ RMAA3188; O.Z. supported by Leverhulme Trust Research Fellowship.
%


%
\end{document}